\documentclass[12pt]{amsart}
\usepackage{geometry}                
\usepackage{graphicx}
\usepackage{amssymb}

\newtheorem{theorem}{Theorem}[section]

\newtheorem{corollary}[theorem]{Corollary}

\theoremstyle{definition}

\newtheorem{question}[theorem]{Question}

\theoremstyle{remark}
\newtheorem{remark}[theorem]{Remark}

\numberwithin{equation}{section}

\newcommand{\cA}{\mathcal{A}}

\newcommand{\D}{\mathbb{D}}

\newcommand{\C}{\mathbb{C}}

\newcommand{\Z}{\mathbb{Z}}

\newcommand{\R}{\mathbb{R}}

\newcommand{\rea}{\operatorname{Re}}

\newcommand{\length}{\operatorname{Length}}

\newcommand{\Arg}{\operatorname{Arg}}

\newcommand{\bd}[1]{\partial #1}

\newcommand{\Mod}{\operatorname{Mod}}

\title{Mapping properties of analytic functions on the disk}
\author{Pietro Poggi-Corradini}
\address{Department of Mathematics, Cardwell Hall, Kansas State University,
Manhattan, KS 66506}
\email{pietro@math.ksu.edu}
\date{}                                           
\begin{document}
\begin{abstract}
There is a universal constant $0<r_0<1$ with the following
property. Suppose that $f$ is an analytic function on the unit disk
$\D$, and suppose that there exists a constant $M>0$ so that the
Euclidean area, counting multiplicity, of the portion of $f(\D)$ which
lies over the disk $D(f(0),M)$, centered at $f(0)$ and of radius $M$,
is strictly less than the area of $D(f(0),M)$. Then $f$ must send
$r_0\overline{\D}$ into $D(f(0),M)$. This answers a conjecture of Don
Marshall.
\end{abstract}

\maketitle

\baselineskip=18pt 

\section{Introduction}\label{sec:intro}
Let $f$ be analytic on $\D$,
$f(0)=0$. For $M>0$, let 
\[\Omega(M)= \{z\in \D:|f(z)|<M\}=\D\cap f^{-1}(D(0,M)).\]
Write $A(R)$ for
the Euclidean area of $f(\Omega(M))$ counting multiplicity, i.e.,
\[A(M)=\int_{\Omega(M)}|f^\prime(z)|^2dA(z).\]
\begin{theorem}\label{cl:main}
There is a universal constant $r_0>0$ such that whenever $f$ is analytic on
$\D$ with $f(0)=0$, and whenever $M>0$ is such that $A(M)<\pi M^2$, then
$f(r_0\D)\subset D(0,M)$.
\end{theorem}  
Theorem \ref{cl:main} was conjectured in \cite{marshall:1989arkiv}
p.~135. As Marshall notes, if $f$ is one-to-one, then two applications
of the $1/4$-Koebe Theorem suffice. In fact, suppose $f$ is conformal
on $\D$ with $f(0)=0$ and $M$ is such that $A(M)<\pi M^2$. Then,
$D(0,M)\setminus f(\D)$ must be non-empty. So if $R>0$ is the largest
radius for which $D(0,R)\subset f(\D)$, then $R<M$. By Koebe, $R\geq
|f^\prime(0)|/4$. Now consider $f^{-1}$ restricted to $D(0,R)$. By
Koebe again, $f^{-1}(D(0,R))$ must contain a disk centered at the
origin of radius $R|(f^{-1})^\prime(0)|/4\geq 1/16$, i.e.,
$f(\frac{1}{16}\D)\subset D(0,R)\subset D(0,M)$. Therefore Theorem
\ref{cl:main} can be thought as a generalization to many-to-one maps
of "two applications of Koebe".

Before proving Theorem \ref{cl:main}, we will reformulate it in a
couple of different ways.  For $0<r<1$, consider the growth function
\[\cA(r)=\frac{A(M(r))}{\pi M(r)^2},\] where
$M(r)=\max_{|z|=r}|f(z)|$. So $\cA(r)$ is a sort of Euclidean Mean
Covering. It is larger than the usual (Euclidean) Ahlfors growth
function.  Then Theorem \ref{cl:main} is equivalent to the following
statements:
\begin{corollary}\label{cor:mean}
There is a universal radius $r_0>0$ so that for every $f$ analytic in
$\D$ with $f(0)=0$, and for every $0<r<r_0$,
\[
\cA(r)\geq 1.
\] 
\end{corollary}
\begin{corollary}\label{cor:mean2}
There is a universal radius $r_0>0$ so that for every $f$ analytic in
$\D$ with $f(0)=0$, and for every $0\leq M\leq \max_{|z|=r_0}|f(z)|$,
\[
A(M)\geq \pi M^2.
\] 
\end{corollary}
\begin{proof}[Proof of Corollary \ref{cor:mean}:]
Let $r_0$ be the constant given in Theorem \ref{cl:main}. Suppose that
for some $r<r_0$, $\cA(r)<1$, i.e. $A(M(r))<\pi M(r)^2$. Then, by
Theorem \ref{cl:main}, $f(r_0\D)\subset D(0,M(r))$. Hence
$M(r_0)\leq M(r)$, but this contradicts the maximum principle.
\end{proof}
Corollary \ref{cor:mean2} is proved similarly.

We end this introduction with a remark and a question related to Theorem \ref{cl:main}.
\begin{remark}
The proof of Theorem \ref{cl:main} provided below breaks down for meromorphic functions when the "Ahlfors island argument" is
used. Indeed, here is a simple counter-example when $f$ is allowed to have poles. The function
\[
f(z)=\frac{z}{\epsilon+2z}
\]
maps $\D$ conformally onto the complement in the extended complex plane of a small disk containing the point $1/2$, $f(0)=0$, and choosing $M=1$ in Theorem \ref{cl:main}, we have $A(1)<\pi$, 
but $f(-\epsilon/3)=-1$, 
for all $\epsilon>0$. 
\end{remark}

\begin{question}
Does a form of Theorem \ref{cl:main} still hold for quasiregular maps
of the unit ball of $\R^n$, $n\geq 2$?
\end{question}
Since the main tool in the proof of Theorem \ref{cl:main} is modulus
of path families, it is safe to expect this phenomenon to persist in
the quasi world. However, a few initial attempts on our part failed
when trying to lift paths.

\noindent {\bf Acknowledgements:} We thank Don Marshall and David
Drasin for helpful comments on a preliminary version of this paper.

\section{A proof of Theorem \ref{cl:main}}\label{sec:pfone}
Without loss of generality $f$ is analytic in $\overline{\D}$
(otherwise consider $f(tz)$ with $t$ close to $1$).  Notice that
$D(0,M)\setminus f(\D)\neq\emptyset$, because otherwise we would have
$A(M)\geq \pi M^2$.  So we can pick a point $\zeta\in D(0,M)\setminus
f(\D)$, and we can rotate and dilate the image so as to have $\zeta>0$
and $M=1$.  We let $G=G(1)$ be the component of $\Omega(1)$ which
contains the origin and we let $r$ be the largest radius so that
$r\D\subset G$.  Then $r\overline{\D}\cap\bd G\neq\emptyset$. Hence
$0\in f(r\overline{\D})\subset \overline{D(0,1)}$
and $f(r\overline{\D})\cap \bd D(0,1)\neq \emptyset$.
Consider a ray $[0,a]$ connecting $0$ to a point $a$ in $r\bd\D\cap\bd G$. Then $E=f([0,a])$ is a curve in $f(r\D)$ connecting $0$ to $\bd D(0,1)$.

Several cases arise. Assume first that $0<\zeta<1/4$, then we consider
the family of circles centered at $\zeta$ with radius $u$ and
$1/4<u<3/4$. Each one of these circles separates $0$ from $\bd
D(0,1)$, so each one must intersect $E$.  Since $E=E(s)=f(sa)$, $0<s<1$ is a piece-wise analytic arc, we can find finitely many open intervals $J_n\subset [0,1]$, $n=1,..,N$, so that $h_n(s)=|E(s)-\zeta|$ is an increasing function  on each interval $J_n$, the ranges $I_n=h_n(J_n)$ are disjoint and their union $\cup_n I_n$ spans the whole interval $[1/4,3/4]$ except for possibly finitely many points. Let $F_n=\{sa: s\in J_n\}$ and $E_n=\{f(sa): s\in J_n\}$. Then
every point  $w_u\in E_n\cap \{|w-\zeta|=u\}$ is hit schlichtly by the corresponding
point $z_u\in F_n$ and a branch of $f^{-1}$ is well-defined in a
neighborhood of $E_n$ such that $f^{-1}(w_u)=z_u$. Also, such branch
can be analytically continued along $\gamma_u(t)=\zeta+u e^{it}$,
for $t>\theta_u$, where $\theta_u=\Arg(w_u)$ plus a multiple of $2\pi$ to be determined later. In fact, $\theta_u$ can be chosen to be a continuous function in $u$ on each interval $I_n$.
We extend
the definition of $\gamma_u(t)$ for $t\in [\theta_u, T_u)$, where
$T_u$ is the largest possible value so as to be able to analytically
continue $f^{-1}$ along $\gamma_u|[\theta_u,T]$ for every $T<T_u$,
with values in $\D$.  Let $\alpha_u(t)$ be the lifted paths
$f^{-1}(\gamma_u(t))$.  We claim that $\alpha_u$ is a simple curve.
If not, there are $t_0\neq t_1$ such that
\[f^{-1}(\gamma_u(t_0))=f^{-1}(\gamma_u(t_1)).\]
Moreover, $t_0$ and $t_1$ can be chosen so that the curve
$f^{-1}(\gamma_u(t))$ between $t_0$ and $t_1$ is a simple closed
curve, and hence bounds a region, which is relatively compact in $\D$
(an "island" in Ahlfors terminology). The image of such region has
boundary in the circle $\{w: |w-\zeta|=u\}$ and hence must wind at least once around some
and therefore all points inside the circle. So $\zeta$ must be in the
image of $f$, by the argument principle, but this is not the case. So
$\alpha_u$ is a simple curve.  As a corollary we find that
$|\alpha_u(t)|$ must tend to $1$ as $t$ tends to $T_u$.

We will estimate the {\sf modulus} $\Mod(\Gamma)$ of the family $\Gamma$ of
all such lifted paths $\alpha_u$. Recall that $\Mod(\Gamma)$ is the
infimum of $\int_{\C}\rho^2(z)dA(z)$ over all {\sf admissible} $\rho$'s,
i.e., positive Borel measurable functions $\rho $ on $\C$ with the
property that
\[\int_{\alpha}\rho ds\geq 1\]
for every path $\alpha\in \Gamma$.

Since $\Gamma$ is a subfamily of all paths connecting $r\overline{\D}$
to $\bd\D$, by the monotonicity of modulus, we get
\[
\Mod(\Gamma)\leq\frac{2\pi}{\log(1/r)}.
\]

On the other hand, $g(z)=\log(f(z)-\zeta)$ can be defined to be an
analytic function on $\D$. So by \cite{ohtsuka1970} \S 2.4 or
\cite{rickman1993} Theorem 8.1 (Poletskii's Inequality),
\[
\Mod(\Gamma)\geq \Mod (g\Gamma)
\]
where $g\Gamma$ is the family of paths of the form $g(\alpha_u)$. In
particular, each path in $g\Gamma$ is a vertical segment of the form
\[\beta_u(t)=\log u+it;  1/4<u<3/4, u\not\in Z, \theta_u+2k\pi <t<T_u+2k\pi\]
for some $k=k(u)\in \Z$. Notice that $k(u)$ can be chosen to be constant for $u$ in an interval $I_n$, and then the set 
$U_n=\cup_{u\in I_n} \beta_u$
is a simply connected  open set on which a branch of $g^{-1}$ is well-defined. In particular, the set $V_n=\cup_{u\in I_n}\alpha_u$ is also open and simply-connected in $\D$ and equals $g^{-1}(U_n)$.

We will obtain a lower bound for $\Mod(g\Gamma)$.
Assume that $\rho$ is a positive Borel function on $\C$
which is admissible for $g\Gamma$. Then
\[
1\leq\int_{\beta_u}\rho dt\leq
\left(\length(\beta_u)\right)^{1/2}\left(\int_{\beta_u}\rho^2
dt\right)^{1/2}
\]

Now, since $U_n$ is an open set, the function $u\mapsto \length(\beta_u)$ is Borel measurable on each interval $I_n$. So changing variables via the coarea formula (see \cite{maly-swanson-ziemer:2003tams} for the coarea formula for Sobolev maps):
\begin{eqnarray}\label{eq:area}
\int_{1/4}^{3/4}\length(\beta_u)udu &  = & \int_{1/4}^{3/4}\int_{\alpha_u} |f(z)-\zeta||g^\prime(z)|dsdu\\ \nonumber
& = & \int_{\cup_n V_n}|f^\prime(z)|^2
dA(z)\leq A(1)<\pi.
\end{eqnarray}

So
\begin{eqnarray*}
\int_{\C}\rho^2 dudt & \geq &
\int_{1/4}^{3/4}\int_{\beta_u}\rho^{2}dtdu \\
& \geq & \frac{4}{3}
\int_{1/4}^{3/4}\frac{udu}{\length(\beta_u)} \\
& \geq & \frac{1}{3\int_{1/4}^{3/4}\length(\beta_u)udu}\geq
\frac{1}{12\pi}.
\end{eqnarray*}
 Putting everything
together
\[
\frac{2\pi}{\log 1/r}\geq \frac{1}{12\pi}
\]
So
\[r\geq \exp\left\{-24\pi^2\right\}.   
\] 
This proves the theorem when $0<\zeta<1/4$.

So assume that $1>\zeta\geq 1/4$ and consider the circles $S_u$
centered at $\zeta$ of radius $u$. Note that
\[M:=\max\{u: E\cap S_u\neq\emptyset\}\geq \zeta \geq
m:=\min\{u: E\cap S_u\neq\emptyset\}\] If $M-\zeta>1/8$, we consider
circles centered at $\zeta$ with radius $\zeta<u<\zeta+1/8$.  Then $E$
intersects every such circle, and we can repeat the same argument as
above. The only difference arises for the estimate (\ref{eq:area}),
because now only a portion of $\gamma_u$ lies schlichtly above $D(0,1)$. But the total length of $\gamma_u$ is comparable to the
length of the part of $\gamma_u$ that lies above $D(0,1)$, i.e. there is a
universal constant $C$, related to the angular measure of $S_u\cap D(0,1)$ about $\zeta$, such that
\[\length(\gamma_u)\leq C\length (\gamma_u\cap D(0,1)).\]

Likewise, if $\zeta-m>1/8$, then the same argument works.

So assume that $\zeta-1/8\leq m<M\leq \zeta+1/8$.
Then $E$ intersects $\bd \D$ in a point $a$ with $\rea a>47/128>1/8$.
Therefore $E$ intersect the boundary $\bd R_u$ of each rectangle 
\[R_u=\{u<x<2\zeta-u, -1-u<y<1+u \}\qquad u\in (0,1/8)\]
Now we can repeat the same argument as above, except that the function
$g$ is replaced by $\log(\phi(f(z))-\zeta)$ where $\phi$ is a
quasiconformal map that fixes $\zeta$ and sends the rectangles $R_u$
to circles centered at $\zeta$ of radius $\zeta+u$. The conclusion
follows because the distortion of $\phi$ is controlled by universal
constants. In this case, one really needs to use Poletskii's inequality
instead of Ohtsuka's result.
 
\section{An alternative approach using Beurling's criterion}

After this paper was written, we realized that one can do away with Poletskii's inequality altogether by using Beurling's criterion. For completeness, we decided to include both proofs.

First we recall Beurling's criterion for extremal metrics, see 
\cite{ahlfors1973} Theorem 4-4, page 61.

\begin{theorem}[Beurling's criterion]\label{thm:beurling}
Let $\Gamma$ be a path family in $\C$.
A Borel measurable function $\rho_0:\C\rightarrow [0,+\infty]$ is extremal (i.e. $\int_\C \rho_0^2dA\leq\int_\C\rho^2dA$ for every admissible $\rho$), if $\Gamma$ contains a subfamily $\Gamma_0$ with the following two properties:
\begin{enumerate}
\item for every $\gamma\in \Gamma_0$:
\[
\int_\gamma \rho_0 |dz| = 1;
\] 
\item for every Borel measurable $h:\C\rightarrow [-\infty,+\infty]$ with the property that 
\[
\int_\gamma h |dz|\geq 0
\] 
for all $\gamma\in \Gamma_0$, we must have
\[
\int_\C h\rho_0 dA \geq 0.
\]
\end{enumerate}
\end{theorem}
The proof of Beurling's theorem is a simple application of the Cauchy-Schwarz inequality.

\begin{proof}[An alternative proof of Theorem \ref{cl:main}:]
We use the same notations as in Section \ref{sec:pfone}. Again we assume first that $0<\zeta<1/4$ and define the curves $\alpha_u$ as before.
We want a lower bound for the modulus of the path family $\Gamma$ of all such lifts $\alpha_u$. Recall also that, except for at most finitely many values of $u$ the paths $\alpha_u$ could be grouped so that  $V_n=\cup_{u\in I_n}\alpha_u$ is an open simply connected set in $\D$, $n=1,2,...,N$.

We claim that the following metric is extremal for $\Gamma$:
\[
\rho_0(z)=\frac{|g^\prime (z)|}{\length(\beta_u)}\qquad\mbox{for $z\in V_n\cap\alpha_u$}
\]
and $\rho_0(z)=0$ otherwise. We let $\Gamma_0=\Gamma$ and check Beurling's criterion.
Property (1) of Theorem \ref{thm:beurling} is clear. For property (2), assume that $h$ is a function as above, then, by the coarea formula,
\begin{equation}\nonumber
\int_\C h\rho_0 dA =  \int_{1/4}^{3/4}\frac{\int_{\alpha_u}h(z)|dz|}{u\length(\beta_u)}du\geq 0.
\end{equation}
So, by coarea and (\ref{eq:area}) 
\begin{eqnarray*}
\Mod(\Gamma) & = & \int_\C \rho_0^2 dA\\
& = & \int_{1/4}^{3/4}\frac{\int_{\alpha_u}|g^\prime(z)||dz|}{u\length(\beta_u)^2}du\\
& = & \int_{1/4}^{3/4}\frac{du}{u\length(\beta_u)}\\
& \geq & \frac{16}{9}\frac{1}{\int_{1/4}^{3/4}\length(\beta_u)udu}>\frac{16}{9\pi}\\
\end{eqnarray*}
This proves Theorem \ref{cl:main} when $0<\zeta<1/4$. The other cases are handled similarly. In the last case when $g$ is a quasiregular map, the extremal metric must be defined with $|g^\prime(z)|$ replaced a.e. by 
\[
\ell(z)=\inf_{|v|=1}|Dg(z) v|
\] 
where $Dg$ is the differential matrix of $g$ and $v$ is a $2\times 1$ vector in $\R^2$.

\end{proof}

\def\cprime{$'$}

\end{document}